\documentclass[journal]{IEEEtran} 
\usepackage{co} 
\usepackage[linesnumbered,lined, algoruled]{algorithm2e}
\usepackage{multirow}
\usepackage{tikz}
\usetikzlibrary{arrows}
\usepackage{hyperref}
\graphicspath{{./}}

\def\minim{\hbox{min }}
\def\maxim{\hbox{max }}

\def\st{\hbox{s.t. }}

\def\pst{\phantom{\st}}
\renewcommand{\Re}{\mathop{\rm Re}\nolimits}
\renewcommand{\Im}{\mathop{\rm Im}\nolimits}

\ifbozze
                 
\fi
\usepackage{enumitem}

\begin{document}

\labelfigure{./}

\title{New results on the local linear convergence of ADMM: a joint approach}

\author{%
Tomaso~Erseghe
\thanks{The author is with the Dipartimento di Ingegneria dell'Informazione, Universit\`a di Padova, Via G. Gradenigo 6/B, 35131 Padova, Italy. Contacts: tel: +39 049 827 7656, fax: +39 049 827 7699, mailto: tomaso.erseghe@unipd.it
} %
\thanks{\copyright 2019 IEEE. Personal use of this material is permitted. Permission from IEEE must be obtained for all other uses, including reprinting/republishing this material for advertising or promotional purposes, collecting new collected works for resale or redistribution to servers or lists, or reuse of any copyrighted component of this work in other works.}
} 

\markboth{SUBMITTED TO IEEE TRANSACTIONS ON SIGNAL PROCESSING}{SUBMITTED TO IEEE TRANSACTIONS ON SIGNAL PROCESSING}
\maketitle
\ifbozze\clearpage\fi


\acrodef{ADMM}{alternating direction method of multipliers} 
\acrodef{APP}{auxiliary problem principle} 
\acrodef{AWGN}{additive white Gaussian noise} 
\acrodef{CPLD}{constant positive linear dependence}
\acrodef{CRLB}{Cramer-Rao lower bound} 
\acrodef{DC}{direct current} 
\acrodef{DP}{dynamic programming} 
\acrodef{IPM}{interior point method} 
\acrodef{ITP}{incomplete-tilted-phi} 
\acrodef{KKT}{Karush Kuhn Tucker} 
\acrodef{LDPC}{low density parity check} 
\acrodef{LLR}{log likelihood ratio} 
\acrodef{LMS}{least mean squares} 
\acrodef{MAP}{maximum a posteriori} 
\acrodef{MEMS}{micro electro-mechanical system} 
\acrodef{ML}{maximum likelihood} 
\acrodef{OPF}{optimal power flow} 
\acrodef{PCC}{point of common coupling} 
\acrodef{PDF}{probability density function} 
\acrodef{PEI}{power electronic interface} 
\acrodef{RMSE}{root mean squared error} 
\acrodef{RF}{radio frequency} 
\acrodef{RSSI}{received signal strength indicator} 
\acrodef{SDP}{semi-definite programming}
\acrodef{SNR}{signal to noise ratio} 
\acrodef{SMG}{smart micro grid} 
\acrodef{SVD}{singular value decomposition}
\acrodef{TOA}{time of arrival} 
\acrodef{TOF}{time of flight} 
\acrodef{WSN}{wireless sensor network}

\begin{abstract}
Thanks to its versatility, its simplicity, and its fast convergence, ADMM is among the most widely used approaches for solving a convex problem in distributed form. However, making it running efficiently is an art that requires a fine tuning of system parameters according to the specific application scenario, and which ultimately calls for a thorough understanding of the hidden mechanisms that control the convergence behaviour. In this framework we aim at providing new theoretical insights on the convergence process and specifically on some constituent matrices of ADMM whose eigenstructure provides a close link with the algorithm's convergence speed. One of the key technique that we develop allows to effectively locate the eigenvalues of a (symmetric) matrix product, thus being able to estimate the contraction properties of ADMM. In the comparison with the results available from the literature, we are able to strengthen the precision of our speed estimate thanks to the fact that we are solving a joint problem (i.e., we are identifying the spectral radius of the product of two matrices) in place of two separate problems (the product of two matrix norms).
\end{abstract}


\section{Introduction}

\IEEEPARstart{I}{n} recent years the \ac{ADMM} has received a considerable attention as an effective method to reach consensus among agents, to implement an optimisation algorithm in a distributed way, or, more generally, to solve networked problems by iteratively applying simple optimisation steps (e.g., see \cite{Boyd:10}). \ac{ADMM} is in close relation with, and in many cases it is equivalent to, a wide number of alternative approaches, e.g., Douglas--Rachford splitting, proximal point algorithms, the split Bregman algorithm, that share the same potential, as well as the same techniques for analysis. 

The reference problem in \ac{ADMM} is 
$$
\eqalign{
&\minim f_1(x_1) + f_2(x_2)  \cr
&\st A_1 x_1 = A_2 x_2 + b\;,
}
\e{AD2}
$$
where $f_i$ are proper, lower semicontinuous, and convex functions, $x_i$ are real valued vectors, and $A_i$ are linear operators (i.e., matrices). Solution to \e{AD2} is found by investigating the associated \emph{augmented Lagrangian} function
$$
\eqalign{
L(x,\lambda) 
 & =  f_1(x_1) + f_2(x_2) + \langle\lambda,  A_1 x_1-A_2 x_2 - b\rangle \cr
 & \hspace*{25mm}+ \fract12\|A_1 x_1-A_2 x_2- b\|_{E}^2\;,
}
\e{AD4}
$$
with $\lambda$ the Lagrange multipliers, $\langle\cdot,\cdot\rangle$ an inner product, and $\|u\|_{E}^2 = \langle u, E u\rangle$ a weighted norm that uses a positive definite \emph{augmentation} matrix $ E\succ 0$. For the sake of ease and practical applicability,  the augmentation matrix $E$ is typically chosen a diagonal positive definite matrix, e.g., $E=\epsilon I$ for some augmentation factor $\epsilon>0$. 

In the above context \ac{ADMM} is an iterative algorithm that solves \e{AD2} by looking for saddle points of the augmented Lagrangian \e{AD4}. This is achieved by performing an alternating direction search where each variable is separately updated by keeping the others fixed, that is, 
$$
\eqalign{
 x_1^+ & = \argmin_{x_1} L( x_1, x_2, \lambda)\cr
 x_2^+ & = \argmin_{x_2} L( x_1^+, x_2, \lambda)\cr
 \lambda^+ & =  \lambda +  E(A_1 x_1^+ - A_2 x_2^+ - b)\;,
}
\e{AD6}
$$
where the $^+$ sign denotes the updated variables. Depending on the choice of functions and linear constraints in \e{AD2}, the above setting enables the possibility of solving a wide range of optimisation problems in distributed form. 

Today \ac{ADMM} is recognised to exhibit a linear (hence fast) convergence provided that parameters are wisely chosen. Unfortunately, parameters selection is a very difficult task. As a consequence, much of the recent effort in \ac{ADMM} is devoted to deriving a deeper understanding of its convergence properties, as well as their dependence on system parameters. The literature on these aspects is vaste. The $O(1/n)$ convergence of \ac{ADMM}, with $n$ denoting the iteration number, is addressed in \cite{he20121, goldstein2014fast, he2015non, davis2016convergence} under weak convexity assumptions. Linear convergence under strong convexity and/or Lipschitz continuity (smoothness) is studied in \cite{ghadimi2015optimal, nishihara2015general, deng2016global, davis2017faster, giselsson2017linear, giselsson2017tight, hong2017linear, liang2017local}, where the derivation is mainly carried out by building inequalities in the dual domain, and by exploiting the properties of proximity functions. A more general statement of linear convergence (under weaker conditions) is available from \cite{hong2017linear}, and a statement on convergence in finite steps under a local polyhedral assumption is available in \cite{liang2017local}. Convergence under specific scenarios is also investigated. Multi-block \ac{ADMM} formalisations, where the target function is split in more than two contributions is available from \cite{lin2014convergence, chen2016direct, deng2017parallel, hong2017linear}. The linear convergence in a decentralised consensus optimisation problem, i.e., a specific problem that can be casted in the general \ac{ADMM} setting to build a distributed algorithm, is proved in \cite{shi2014linear, jakovetic2015linear, iutzeler2016explicit}. Convergence in the presence of inexact solvers is studied in \cite{deng2016global}, while a perspective on non-convex problems is available from \cite{attouch2013convergence, wang2015global, li2016douglas}. 

Taken from the perspective of the usability of the result, a notable approach among the considered ones is given by the works of Giselsson and Boyd \cite{giselsson2017linear, giselsson2017tight}, and of Liang, Fadili, and Peyr\'e of \cite{liang2017local}, which we aim at generalising in some key aspects, and, up to a certain extent, also at merging and at simplifying. More specifically, in our work we provide a number of insights that are useful both on a practical as well as on a theoretical level, namely: we show that the Douglas--Ratchford structure is available in \ac{ADMM} directly in the primal domain, and not only in the dual domain; we simplify derivations by avoiding side concepts (e.g., co-coercivity), and only relying on simple properties of the Legendre Fenchel transform, which, incidentally, also ensures a straightforward understanding of the convergence process in the presence of non-smooth functions; we relax \cite[Assumption~2.ii]{giselsson2017linear}, which limits the applicability of the results; we jointly consider both the functions of the alternate direction process and accordingly develop an efficient method to locate the eigenvalues of a (symmetric) matrix product; we use and generalise the findings of \cite{liang2017local}, thus being able to effectively estimate the convergence speed in the vicinity of the target point.  

The paper is organised as follows. \sect{RE} provides the basis of our approach by showing that a link to a (generalised) Douglas--Rachford splitting is available directly in the primal domain, and that its properties are closely linked to a Legendre--Fenchel transform. These provide the essential structure of \ac{ADMM} which is used in \sect{CN} to establish an equivalent to \cite[Corollary~2]{giselsson2017linear} under a more general framework where a full-row-rank matrix assumption is dropped. Insights under the presence of discontinuities in the derivatives of functions $f_i$ are also provided. In \sect{JO} we study the interdependencies between the two alternating steps in \e{AD6}, and identify a specific matrix $R$ as the key feature of \ac{ADMM}. The in-depth study of the eigenvalue structure of $R$, which controls the convergence process, is given in \sect{LO}, while its (long) technical proof is available in the Appendix. A brief example of application is given \sect{EX}. \sect{CO} concludes the paper.

\Section[RE]{Essential structure of \ac{ADMM} }

\subsection{Douglas--Rachford splitting structure in the primal domain}

The starting point is to reveal an alternative representation of \ac{ADMM} in the form of a (\emph{generalised}) Douglas--Rachford splitting, which is particularly well suited for studying the convergence behaviour. The Douglas--Rachford representation is a standard result of the \ac{ADMM} literature, commonly adopted in the dual domain. Here, instead, we reveal it directly in the primal domain. To this aim, in the following we use a compact notation where 
$$
\eqalign{
\C P_i(u) & = \argmin_{x_i}  f_i(x_i) +  \fract12\| A_i x_i - u\|_E^2\;,\cr
\C D_i(u) & = 2E^{\frac12}A_i \C P_i(E^{-\frac12}u) - u + \cases{-E^{\frac12}b& $i=1$\cr +E^{\frac12}b & $i=2$\,,}
}
\e{AD8}
$$
are, respectively, a \emph{proximity} operator and a \emph{reflected proximity} operator for function $f_i$, $i\in\{1,2\}$. With the above notation, the \ac{ADMM} updates \e{AD6} can be written as
$$
\eqalign{
 x_1^+ & = \C P_1(A_2 x_2 + b  -\tilde\lambda)\cr
 x_2^+ & = \C P_2( A_1 x_1^+ - b + \tilde\lambda)\cr
\tilde\lambda^+ & =  \tilde\lambda +  A_1 x_1^+ - A_2 x_2^+ - b\;,
}
\e{AD10}
$$
where $\tilde\lambda=E^{-1}\lambda$. This reveals that the update of $x_i$ is constrained by the \emph{proximity} operators $\C P_i$, which ensures a form of coordination among the alternate minimisation steps. 

Mimicking the most efficient solutions available in the literature, we further generalise \e{AD10} in a so called \ac{ADMM} \emph{with scaled variables} where the updates are modified in the form
$$
\eqalign{
 x_1^+ & = \C P_1(A_2 x_2 + b -\tilde\lambda)\cr
 y^+ & = 2q  A_1 x_1^+ + (1-2q)  (A_2  x_2+b)\cr
 x_2^+ & = \C P_2(y^+ -b + \tilde\lambda)\cr
\tilde \lambda^+ & = \tilde \lambda +  y^+- A_2 x_2^+ - b\;,
}
\e{AD12}
$$
for some \emph{real} constant $q\neq0$. Note that in \e{AD12} the variable $y$ plays the role of an estimate of $A_1x_1$ (or, equivalently, of $A_2x_2+b$), exploiting both $x_1$ and $x_2$ in a (possibly) more efficient way that linearly combines both the directions of the alternate search. The above corresponds to the standard \ac{ADMM} algorithm when $q=\frac12$. 

The Douglas--Rachford structure becomes evident by introducing variable 
$$ 
z^+ =  E^{\frac12}(y^+ -b + \tilde\lambda)
\e{VAR2}
$$ 
in place of variable $\tilde\lambda$. Note that, although the presence of $E^{\frac12}$ could be dropped, this will guarantee a useful form of symmetry in later derivations, and therefore we keep it despite the redundant notation. Hence, by replacing $\tilde\lambda$ with $z^+$, the iterative algorithm \e{AD12} can be expressed in the equivalent form\footnote{Note that: the first of \e{AD10c2} is obtained by simply replacing \e{VAR2} in the third of \e{AD12};  the third of \e{AD10c2} is obtained from the fourth of \e{AD12} by replacing $\tilde\lambda$ according to \e{VAR2}, and then $y^+$ according to the second of \e{AD12}; the second of \e{AD10c2} is obtained from the first of \e{AD12} by replacing $\tilde\lambda$ according to \e{VAR2}, then $z^+$ according to the third of \e{AD10c2}, and finally $y^+$ according to the second of \e{AD12}.}
$$
\eqalign{
 x_2 & = \C P_2 (E^{-\frac12}z)\cr
 x_1^+ & =  \C P_1 (2A_2x_2+b-E^{-\frac12}z)\cr
 z^{+} & =  z + 2q\,E^{\frac12}  ( A_1 x_1^+ - A_2 x_2-b)\;,
}
\e{AD10c2}
$$
showing that $z$ plays the role of a \emph{state} variable in the \ac{ADMM} update. Furthermore, by replacing the updates of $x_1$ and $x_2$ in the third of \e{AD10c2}, then the update of $z$ can be more compactly expressed through the (generalized) Douglas--Rachford splitting relation
$$
z^{+} = \C R(z) = (1-q) \, z +  q \,\C D_1 \C D_2 ( z)\;,
\e{AD10d2}
$$
where $\C D_i$ are the reflected proximity operators defined in \e{AD8}. We also have
$$
\eqalign{
 x_2 & = \C P_2 (E^{-\frac12}z)\cr
 x_1^+ & =  \C P_1 (E^{-\frac12}\C D_2(z))\cr
}
\e{AD10c4}
$$
which express the link between $z$ and the target variables that we are estimating. Note that \e{AD10d2} evidences the essential structure of \ac{ADMM} directly in the primal domain (e.g., note the closeness with its dual counterpart \cite[eq.~(24)]{giselsson2017linear}), and shows that convergence properties are closely linked to the contraction properties of operator $\C D_1 \C D_2$. 

Formally, the three approaches \e{AD10}, \e{AD12}, and \e{AD10c2} share the same stationary points, in the sense that there is a one-to-one correspondence between their stationary points. The relation is given by the mappings $y^*=A_1x_1^*=A_2x_2^*+b$ and $z^*=E^{\frac12}(y^*-b+E^{-1}\lambda^*)$, with $^*$ denoting the stationary values. In addition, these stationary points coincide with those of the Lagrangian \e{AD4}, i.e., they solve the minimisation problem. Evidently, we are implicitly assuming that such stationary points (i.e., a solution) exist, that is throughout this paper we are taking the following
\begin{assumption}\label{ass-exist}
A solution to \e{AD2} exists, i.e., equation \e{AD10d2} has at least one stationary point.\hfill~$\Box$
\end{assumption}

\subsection{Legendre--Fenchel structure of $\C D_i$}

According to \e{AD10d2}, the convergence behaviour of the \ac{ADMM} algorithm can be understood (and optimised) by investigating the properties of operators $\C D_i$. In this respect, it is appropriate to give $\C D_i$ a neater form, given by a Legendre--Fenchel transform, which better reveals the dependence on (the properties of) functions $f_i$. The techniques used in the following are taken from the state-of-the-art literature on \ac{ADMM}, and wisely exploited.

The starting point for inspecting operator $\C D_i$ is to first interpret the structure of the proximal operator $\C P_i$. To do so we introduce the shorthand notations
$$
F_i = E^{\frac12}A_i\;,\quad M_i = F_i^T F_i = A_i^T E A_i\;,
\e{VAR4}
$$
and 
$$
\eqalign{
q_M(u)  & = \fract12\|u\|_M^2\cr
g_i (x_i) & = f_i (x_i)+ q_{M_i}(x_i)\;,\cr
}
\e{VAR9}
$$
which allow writing
$$
\eqalign{
\C P_i(E^{-\frac12}u) & = \argmin_{x_i} f_i(x_i) + \fract12\|A_i x_i-  E^{-\frac12}u \|^2_E\cr
 & =  \argmin_{x_i}  g_i (x_i) - \langle  x_i,  F_i^Tu  \rangle\;.
}
\e{VAR6}
$$
In oder to build upon \e{VAR6}, an important requirement is given by
\begin{assumption}\label{ass-convex} Functions $g_i$ are proper, lower semicontinuous, and convex.\hfill$\Box$
\end{assumption}
This is naturally satisfied if, as discussed in the Introduction, functions $f_i$ are proper, lower semicontinuous, and convex. However, Assumption~\ref{ass-convex} is more general, and further considers the possibility that $f_i$ are non-convex, provided that the structure of matrices $M_i$ is such to make the sum $f_i + q_{M_i}$ convex. 

We then observe that the second line of \e{VAR6} is related to the Legendre--Fenchel transform $g_i^*$ of $g_i$, defined as (e.g., see \cite[\S11]{rockafellar2009variational}) 
$$
g_i^*(u) = \sup_{x_i}\, \langle  x_i,  u  \rangle - g_i(x_i)\;.
\e{VAR10}
$$
The explicit link between \e{VAR10} and \e{VAR6} is the standard result \cite[Proposition~11.3]{rockafellar2009variational}
$$
\partial g_i^* (u) = \left\{ \argmax_{x_i}\, \langle  x_i,  u  \rangle - g_i(x_i) \right\}\;,
\e{VAR12}
$$
where $\partial g_i^* $ expresses the sub-gradient of $g_i^*$, and the curly brackets identify the sets of arguments that maximise their target function. The validity of \e{VAR12} is guaranteed by Assumption~\ref{ass-convex}. Therefore, the explicit link between $g_i^*$ and the proximity operator is, by inspection,
$$
\C P_i(E^{-\frac12}u) \in \partial g_i^* (F_i^Tu)\;,
\e{VAR14a0}
$$
which further provides
$$
\C D_i(u) \in 2F_i \partial g_i^* (F_i^Tu) - u \mp E^{\frac12}b\;.
\e{VAR14}
$$
Equation \e{VAR14} shows that the contraction properties of operators $\C D_i$ are (linearly) linked to those of the Legendre--Fenchel transform $g_i^*$. Hence, an understanding of the convergence behaviour of \ac{ADMM} needs investigating $\partial g_i^*$.

\Section[CN]{Contraction properties of the reflected proximity operators $\C D_i$}

In the following we discuss the properties of \e{VAR14} by first concentrating on a simplified context where functions $g_i$ are convex and smooth, which essentially leads to the outcomes of \cite{giselsson2017linear}, although under a more general framework. We then generalise the result to a more realistic context where functions are non-smooth, e.g., in the presence of piecewise functions or boundaries, which clarifies the added value of the proposed approach.

\subsection{Smooth functions scenario}

We first take the simple (but effective) scenario used in \cite{shi2014linear, giselsson2017linear}, namely 

\begin{assumption}\label{ass-smooth} 
Function $f_i$ is $S_i$-smooth and $C_i$-strongly convex, meaning that both functions $f_i-q_{C_i}$ and $q_{S_i}-f_i$ are convex, where $S_i$ and $C_i$ are symmetric matrices satisfying $S_i\succeq C_i$, and where $\succeq$ is the positive semidefinite ordering operator.\hfill$\Box$
\end{assumption}

According to Assumption~\ref{ass-convex}, this corresponds to a Lipschitz property that limits the Hessian of $g_i$ to the range (see also the definition of $g_i$ in \e{VAR9})
$$
0 \preceq C_i + M_i \preceq \partial^2 g_i \preceq S_i  + M_i\;.
\e{VAR20}
$$
As a consequence of \cite[Theorem~18.15]{bauschke2011convex}, the above reveals a Lipschitz property in the Legendre--Fenchel domain, namely the property 
$$
0 \preceq (S_i+M_i)^{-1} \preceq \partial^2 g_i^* \preceq (C_i+M_i)^{-1}\;,
\e{VAR24}
$$
where the negative exponent stands for a pseudo--inverse operator when matrices are non--invertible. This, in connection with \e{VAR14}, evidences that the sub-gradient of operator $\C D_i$ satisfies the Lipschitz property
$$
-I \preceq L_i  \preceq \partial D_i \preceq U_i\;,
\e{VAR26}
$$
where $I$ is the identity matrix, and where
$$
\eqalign{
L_i & = 2F_i(S_i+M_i)^{-1} F_i^T - I\cr
U_i & = 2F_i (C_i+M_i)^{-1}F_i^T - I 
}
\e{VAR28}
$$
are lower and upper bounds. Incidentally, the above (as well as what follows) also holds in the limit for matrices $C_i+M_i$ and $S_i+M_i$ having eigenvalues which are arbitrarily small (i.e., $\rightarrow0$) or arbitrarily large (i.e., $\rightarrow \infty$). 

A neater expression for bounds \e{VAR28} can be obtained by introducing matrices
$$
\tilde{F}_i = F_i M_i^{-\frac12}\;,
\e{UU2}
$$
providing\footnote{This is a straightforward consequence of standard properties of the Moore--Penrose pseudo--inverse $F_i^\dag=M_i^{-1}F_i^T$.}
$$
\eqalign{
L_i & = \tilde{F}_i  \left( \frac{I-\tilde{S}_i}{I+\tilde{S}_i} \right) \tilde{F}_i^T - (I-\tilde{F}_i\tilde{F}_i^T)\cr
U_i & =  \tilde{F}_i  \left( \frac{I-\tilde{C}_i}{I+\tilde{C}_i} \right) \tilde{F}_i^T - (I-\tilde{F}_i\tilde{F}_i^T) \;,
}
\e{VAR28b}
$$
with symmetric matrices
$$
\eqalign{
\tilde{S}_i &= M_i^{-\frac12} S_i M_i^{-\frac12}\cr
\tilde{C}_i &= M_i^{-\frac12} C_i M_i^{-\frac12}
}
\e{VAR28c}
$$
satisfying $\tilde{S}_i\succeq\tilde{C}_i \succeq -I$ (see also \e{VAR20}). Note that $L_i$, $U_i$, $\tilde{S}_i$, and $\tilde{C}_i$ are symmetric real matrices that are $\succeq-I$, and their eigenvalues are, therefore, real valued and greater than or equal to $-1$. 

The structure revealed in \e{VAR28b} is common to a number of findings available in the literature \cite{lin2014convergence, deng2016global, giselsson2017linear}. However, a careful reader can appreciate that we leverage the same essential properties governing \ac{ADMM}, but the present derivation avoids the difficulties involved with a study performed in the dual domain, and simplifies the proofs. In addition, the present formalisation contains some generalisations that can be fully appreciated by investigating the eigen-structure of matrices. Specifically, from standard linear algebra considerations it is easy to see that the eigenvalues of lower and upper bound matrices \e{VAR28b} are given, respectively, by
$$
\eqalign{
\ell_{i,k} & = \cases{h(\tilde{s}_{i,k}) & \hbox{$k=1,\ldots,{\rm rank}(A_i)$}\cr
	\rule{0mm}{5mm}-1 & \hbox{otherwise,}}\cr
\rule{0mm}{9mm}\nu_{i,k} & = \cases{h(\tilde{c}_{i,k}) & \hbox{$k=1,\ldots,{\rm rank}(A_i)$}\cr
	\rule{0mm}{5mm}-1 & \hbox{otherwise,}}\cr
}
\e{HH10}
$$
where (see a graphical representation in \fig{CO2}\Fig[t]{CO2})
$$
h(x) = \frac{1-x}{1+x}\;,
\e{HH11}
$$
and where $\tilde{s}_{i,k}$ and $\tilde{c}_{i,k}$ are those (real valued) eigenvalues of matrices $\tilde{S}_i$ and $\tilde{C}_i$, respectively, that are associated with the span of matrix $M_i$, and are therefore in number of ${\rm rank}(M_i)={\rm rank}(A_i)$. Incidentally, the eigenvalues $\tilde{s}_{i,k}$ and $\tilde{c}_{i,k}$ can be extracted, respectively, from matrices $S_iM_i^{-1}$ and $C_iM_i^{-1}$ which avoid the square--root operation and are diagonalisable by construction. All the eigenvalues in \e{HH10} are real valued. 

It is important to observe that the eigenvalues $-1$ in \e{HH10} are a consequence of the presence in \e{VAR28b} of the contribution $I-\tilde{F}_i\tilde{F}_i^T$ (if active), which is the projector associated to the kernel of $F_i$. These values are the most significant novelty we are providing, and are also a straightforward consequence of the fact that we are avoiding \cite[Assumption~2.ii]{giselsson2017linear}, i.e., the request of $A _i$ being a full-row-rank matrix. Incidentally, the presence of eigenvalues $-1$ is common in the very many optimisation problems where some of the variables are duplicated more than once.

\subsection{Smoothness + convexity scenario}

If we focus on a scenario where functions $f_i$ are convex, which corresponds to setting $C_i\succeq0$, then the natural consequence is that eigenvalues $\tilde{s}_{i,k}$ and $\tilde{c}_{i,k}$ are non negative. Hence, by inspection of the structural function $h(x)$ for $x\ge0$ (see \fig{CO2}), all eigenvalues $\ell_{i,k}$ and $\nu_{i,k}$ as defined in \e{HH10} are constrained in the range $[-1,1]$. The straightforward outcome is that 
$$
U_i \preceq I
\e{VAR30}
$$
holds, which further guarantees that $-I \preceq \partial D_i \preceq I$, and, in turn, ensures the continuity of functionals $\C D_i$ as well as their \emph{non expansiveness}. 

If we further denote with $\ell_{i,\min}\ge-1$ the minimum eigenvalue of $L_i$, and with $\nu_{i,\max}\le1$ the maximum eigenvalue of $U_i$, we also have $I \cdot \ell_{i,\min}  \preceq \partial D_i \preceq I \cdot \nu_{i,\max}$, whose specific consequence is the contractive property
$$
\|\C D_i(u)-\C D_i(v)\| \le  \mu_i \|u-v\|\;,
\e{VAR32}
$$
with contraction parameter 
$$
\eqalign{
\mu_i & = \max(|\ell_{i,\min}|,|\nu_{i,\max}|)\cr
 & = \max(-\ell_{i,\min},\nu_{i,\max})
 }
\e{VAR34}
$$
satisfying $0<\mu_i\le 1$. Under the assumption of a full-row-ranked matrix $A_i$, and although we are using a different notation, this result perfectly corresponds to \cite[Proposition~5]{giselsson2017linear} (the more general result of the paper).

\subsection{Non-smooth functions scenario}

In very many applications the functions $f_i$ involved in \ac{ADMM} do not satisfy the smoothness property, either because they are non smooth by construction, or because they have a limited domain. In these situations, the formalisation identified by \e{VAR14} is able to fully (and simply) capture the essence of the problem by a simple calculation of a Legendre-Fenchel transform. To make our point, we investigate a context where variables are simply duplicated, and functions $f_i$ are expressed as a linear combination of one--dimensional \emph{convex} contributions. This is equivalent to considering that $E=\epsilon$ and $A_i=a$ are scalar values, that $b=0$ (for simplicity, since the value of $b$ is irrelevant for the contracting properties of $\C D_i$), and that functions $f_i(x)$ are one-dimensional. This simplified scenario is able to capture the essence of the most general case.

To begin with we assume a piecewise linear form for $f_i$, thus obtaining a tilted-staircase reflected proximity operator $\C D_i$ as illustrated in \fig{CO4}\Fig[t]{CO4}. The expression/structure of $\C D_i$ was derived by application of \e{VAR14}, and by exploiting the identity $\partial g_i^* = (\partial g_i)^{-1}$, which is valid thanks to the convexity of $g_i$ \cite[Proposition~11.3]{rockafellar2009variational}. The relevant facts that can be observed in \fig{CO4} are:
\begin{enumerate}[label=\alph*)]
\item the operator $\C D_i$ has a (rotated) staircase structure where a negative $-1$ slope alternates with a positive $+1$ slope, negative slopes being associated to breaking points $x_k$ and positive slopes to the constant derivatives $m_j$, the general expression of lines being available in \fig{CO4}; observe that in a multidimensional (generalised) context, this would map into the fact that the domain is partitioned in regions where slopes can either be $-1$ or $+1$, where $1$ is intended as the all-ones vector;
\item the larger the augmentation variable $\epsilon$ the wider the positive--sloped segments, while the smaller $\epsilon$ the wider the negative--sloped segments, a result which is a consequence of the value 
$$
v_{k,j} = \sqrt{\epsilon}ax_k+\frac{m_j}{\sqrt{\epsilon}a}
\e{DD4}
$$ 
of the coordinate where a negative--sloped ($x_k$) and a positive--sloped ($m_j$) curve intersect (see the exemplifying $v_{4,4}$ and $v_{5,4}$ in \fig{CO4});
\item the active range $[\ell_{i,\min},\nu_{i,\max}]\subset[-1,1]$ (see \sect{CN}.B) may significantly change over different intervals, e.g., in the depicted case where in the region around $v$ we have $\nu_{i,\max}=1$ but $\ell_{i,\min}>-1$; as a natural consequence, the convergence behaviour strongly depends on the limit point.
\end{enumerate}
The above provides an operative method to study convergence bounds $[\ell_{i,\min},\nu_{i,\max}]$ in the vicinity of a target point, e.g., in the vicinity of the target $z^*$, when this is known or it can be calculated with some precision.

The findings of \fig{CO4} can be generalised to build a simple procedure that identifies the proximity operator $\C D_i$ when functions $f_i$ are non-smooth, but not necessarily piecewise-linear. To do so, observe that, according to \fig{CO4}, the value $v_{k,j}$ in \e{DD4} is mapped by the reflected proximity operator into
$$
\C D_i(v_{k,j}) = \sqrt{\epsilon}ax_k-\frac{m_j}{\sqrt{\epsilon}a}\;.
\e{DD4b}
$$
If a piecewise--linear curve like the one in \fig{CO4} is taken to approximate, in the limit, any given function $f_i$, then, by virtue of continuity, we can interpret value $m_j$ in \e{DD4} and \e{DD4b} as $m_j=f'_i(x_j)$ and we are  therefore able to capture the shape of function $\C D_i(u)$ through the couple of functions
$$
\eqalign{
u(x) & = \sqrt{\epsilon}ax+\frac{f'_i(x)}{\sqrt{\epsilon}a}\cr
\C D_i(u(x)) & = \sqrt{\epsilon}ax-\frac{f'_i(x)}{\sqrt{\epsilon}a}
}
\e{DD6}
$$
as well as its derivative through
$$
\partial\C D_i(u(x))  = \frac{\partial \C D_i(u(x))/\partial x}{\partial u(x)/\partial x} = h\left(  \frac{f''_i(x)}{\epsilon a^2} \right)\;.
\e{DD8}
$$
This identifies a simple and effective method to practically infer the structure of $\C D_i$ when $f_i$ is a piecewise function, where:
\begin{enumerate}
\item \e{DD6} and \e{DD8} provide the shape of $\C D_i$ associated with the smooth pieces of a piecewise function $f_i$;
\item the link between such pieces is captured, for each breakpoint $x_k$, by a negative--sloped contribution of the form $2\sqrt{\epsilon}ax_k-u$ (see \fig{CO4}).
\end{enumerate}
Although derived for a convex function, the above method is perfectly applicable to a non--convex setting, provided that Assumption~\ref{ass-convex} holds, i.e., that $u(x)$ is an increasing function.

\Section[JO]{Joint role of the reflected proximity operators $\C D_1\C D_2$}

\subsection{Merging the reflected proximity operators}

The fact that, in the reference equation \e{AD10d2}, the reflected proximity operators appear in the form $\C D_1\C D_2$ suggests an analysis of their joint role. To do so we assume that matrices $L_i$ and $U_i$ commute, which allows writing the bounds in \e{VAR26} in the form
$$
\eqalign{
L_i & = V_i^T \mathrm{diag}(\ell_i) V_i\cr
U_i & = V_i^T \mathrm{diag}(\nu_i) V_i\;.
}
\e{VAR36}
$$
for some \emph{orthogonal} (real--valued and unitary) matrix $V_i$, where $\ell_i$ and $\nu_i$ collect the eigenvalues \e{HH10}, and where, thanks to assumption $L_i\preceq U_i$, the eigenvalues satisfy the ordering
$$
\ell_{i,k}\le\nu_{i,k}\;.
\e{VAR36b}
$$ 
We are therefore taking the following:
\begin{assumption}\label{ass-commute}
Matrices $\tilde{S}_i$ and $\tilde{C}_i$ in \e{VAR28c} commute, or, equivalently, $S_iM_i^{-1}C_i=C_iM_i^{-1}S_i$.\hfill~$\Box$
\end{assumption}
\paragraph*{Remark 1} The assumption is satisfied in very many practical cases. Furthermore, it can be always ensured by simply relaxing the true (stricter) bounds declared by $C_i$ and $S_i$. In practical setups this is straightforwardly obtained, e.g., by making sure by construction that $M_i$, $C_i$, and $S_i$ are diagonal matrices.\hfill~$\Box$
\paragraph*{Remark 2} In the most general case where $L_i$ and $U_i$ do not commute, \e{VAR36} holds in the form 
$$
\eqalign{
L_i & = V_i^T \mathrm{diag}(\breve{\ell}_i) V_i\cr
U_i & = V_i^T \mathrm{diag}(\breve{\nu}_i) V_i\;,
}
\e{MLL2}
$$
where $V_i$ is real--valued and invertible, but not necessarily unitary, and where $\breve{\ell}_i$ and $\breve{\nu_i}$ are vectors identifying the common diagonal structure between matrices $L_i$ and $U_i$. This can be proved, e.g., by exploiting the fact that two positive--semidefinite matrices can be simultaneously diagonalised by an invertible matrix (e.g., see \cite{newcomb1961simultaneous}). Thanks to the Inertia Theorem, $\breve{\ell}_i$ and $\breve{u}_i$ keep the same signs of $\ell_i$ and $u_i$, and some relations on ranges can also be inferred. However, the fact that $V_i$ is non--unitary makes some of the derivations of this paper hard, this being the motivation for taking Assumption~\ref{ass-commute}. \hfill~$\Box$

The above allows to characterise the Lipschitz properties of $\C D_i(u)$. Specifically, it is easy to verify that the operator 
$$
\tilde{\C D}_i(u)=V_i \C D_i(V_i^T u)
\e{VAR38a}
$$ 
satisfies the Lipschitz property\footnote{This is a consequence of the fact that $\partial\tilde{\C D}_i=V_i \partial\C D_i  V_i^T$.} $\mathrm{diag}({\ell}_i) \preceq \partial \tilde{\C D}_i \preceq \mathrm{diag}({\nu}_i)$. This corresponds to verifying that the relation
$$
\tilde{\C D}_i(u) = \tilde{\C D}_i(v) + \mathrm{diag}(\alpha_i) (u-v)
\e{VAR44}
$$
holds for some vector $\alpha_i\in[{\ell}_i,{\nu}_i]$, which is a straightforward consequence of observing (along each direction) an upper and a lower bound on the derivatives. By mapping back onto $\C D_i$ we further have
$$
\C D_i(u) -\C D_i(v) =  V_i^T\,\mathrm{diag}(\alpha_i) V_i (u-v)\;,
\e{VAR46}
$$
for the same vector $\alpha_i$. The use of \e{VAR46} in \e{AD10d2} finally provides our target result
$$
\eqalign{
& \C R(u)-\C R(v)   =  R(\alpha)  \cdot  (u-v)\cr
& \rule{0mm}{5mm}\hspace*{10mm} R(\alpha)  = (1-q) \, I+  q\, N(\alpha)\cr
& \hspace*{10mm} N(\alpha)  = V_1^T\,\mathrm{diag}(\alpha_1) \,V_1\, V_2^T\,\mathrm{diag}(\alpha_2) \,V_2\;,
}
\e{VAR48}
$$
for some $\alpha=[\alpha_1,\alpha_2]$ and $\alpha_i\in[{\ell}_i,{\nu}_i]$. Note that, although expressed with a different notation, the structure of $R(\alpha)$ in \e{VAR48} corresponds to \cite[Eq.~(16)]{liang2017local}, with the added value that $\alpha$ is not necessarily capturing the behaviour only in the strict vicinity of the target point $z^*$, but rather on a larger span. This corresponds to the fact that $\alpha$ belongs to a range as opposed to taking a specific value.

\subsection{Global contraction parameter}

The above analysis suggests that the contraction parameter of global interest, $\mu$, can be identified as the solution to the non-convex problem
$$
\mu =   \max_{\alpha_i\in[{\ell}_i, {\nu}_i]} \,\|R(\alpha)\|\;,
\e{VAR52}
$$
where the matrix norm $\|\cdot\|$ identifies the largest singular value of $R(\alpha)$. By applying this result to \e{VAR48} we obtain
$$
\|\C R(u)-\C R(v)\|   \le \mu   \| u-v\|\;,
\e{VAR48ter}
$$
so that, by respectively choosing $(u,v)=(z,z^-)$ and $(u,v)=(z,z^*)$, with $z^*$ the target value, this establishes that 
$$
\eqalign{
\|z^+-z\| & \le \mu \|z-z^-\|\cr 
\|z^+-z^*\| & \le \mu \|z-z^*\|\;,
}
\e{VAR50}
$$
which corresponds to a contractive property of the \ac{ADMM} update when $\mu<1$. 

Note that this generalises the analysis of \cite{giselsson2017linear} where contraction properties \e{VAR32}-\e{VAR34} are independently exploited to reveal that
$$
\|\C R(u)-\C R(v)\| \le  \,\underbrace{\big((1-q) + q \mu_1\mu_2\big)}_{\tilde{\mu}}\cdot \|u-v\|\;,
\e{VAR32gen}
$$
and it is easy to verify that $\tilde{\mu}\ge\mu$ thanks to the fact that we are jointly considering the role of the reflected proximity operators. We are also generalising (to \ac{ADMM}) the results available for a Douglas--Rachford splitting approach in \cite{giselsson2017tight}, and in fact we are able to account for the presence of matrices $A_i$, as well as to deal with a non-elementary matrix structure for $S_i$ and $C_i$. The following example shows the quality of the proposed methodology through an improvement over available results, even in a very simple case.

\paragraph*{Example} We consider the problem of \cite[Sect.~5]{giselsson2017tight}, that is a convex scenario with $C_1=0$, $S_1 = \beta I\succeq0$, $C_2= \sigma I\succeq0$, $S_2=\infty I$, $A_1=A_2=I$, and $E=\gamma I$. Use in \e{VAR28b} provides
$$
\eqalign{
L_1 & = h(\fract\beta\gamma)\,I\;,\cr
U_1 & = I\;,\cr
}\quad \eqalign{
L_2 & = -I\;,\cr
U_2 & = h(\fract\sigma\gamma)\,I\;,\cr
}
\e{JJ30a}
$$
so that $V_1=V_2=I$, and in \e{VAR48} we have $N(\alpha)=\alpha_1\alpha_2 I$ with $h(\fract\beta\gamma)\le\alpha_1\le1$ and $-1\le\alpha_2\le h(\fract\sigma\gamma)$. As a consequence, the contraction parameter \e{VAR52} becomes
$$
\eqalign{
\mu =\; &\!  \max_{\alpha} \,\big|(1-q) + q \, \alpha_1\alpha_2\big| \cr
&\st \alpha_1\in[h(\fract\beta\gamma), 1]\;,\; \alpha_2\in[-1,h(\fract\sigma\gamma)]\;.
}
\e{JJ30}
$$
One can easily see that the effective range of the product $\alpha_1\alpha_2$ is captured by the expression
$$
-1 \le \alpha_1\alpha_2 \le H_1 = \max\Big(	h(\fract\sigma\gamma), 
	\;h(\fract\beta\gamma)h(\fract\sigma\gamma),\; -\!h(\fract\beta\gamma) \Big)\;,
$$
which can be used to identify the optimum $\gamma$, i.e., the one providing the smallest range. With a little effort\footnote{We skip the derivation since it is lengthy but trivial. As a matter of fact it simply implies testing all the possible signs of $h(\fract\beta\gamma)$ and $h(\fract\sigma\gamma)$, along with the constraints $\sigma<\beta$ or $\sigma>\beta$.} we find that this optimal parameter is $\gamma^*=\sqrt{\sigma\beta}$ providing 
$$
H_1^* = \cases{
\displaystyle\frac{2}{1+\sqrt{\frac\sigma\beta}} - 1 \ge 0& $,\;\sigma<\beta$\cr
\displaystyle\frac{2}{1+\frac12\sqrt{\frac\sigma\beta}+\frac12\sqrt{\frac\beta\sigma}} - 1 \le 0& $,\;\sigma>\beta$.
}
$$
Use of this result in \e{JJ30} further ensures that the minimum contraction parameter is $\mu^*=2q^*-1$ for $q^*=2/(3-H_1^*)$. Hence, by substitution, the optimal contraction parameter in the considered setting is
$$
\mu^* =  \cases{
\displaystyle\left(1 +2\sqrt{\fract\sigma\beta}  \right)^{-1} & $,\;\sigma<\beta$\cr
\displaystyle\left(1  +\sqrt{\fract\sigma\beta} +\sqrt{\fract\beta\sigma} \right)^{-1} & $,\;\sigma>\beta$.
}
\e{JJ32}
$$
We are therefore able to strengthen the result of \cite[Theorem~5.6]{giselsson2017tight}, which identifies the same optimal parameter $\gamma ^*$, but a weaker bound $(1+\sqrt{\sigma/\beta})^{-1}>\mu^*$. The difference is depicted in \fig{CO6}.\Fig[t]{CO6}\hfill~$\Box$

\Section[LO]{Characterisation of the local convergence through eigenvalues}


Following \cite{liang2017local}, an alternative to \e{VAR52} is to investigate the measure
$$
\rho = \max_{\alpha_i\in[{\ell}_i, {\nu}_i]} \,\rho(R(\alpha)) \le \mu \;,
\e{VAR52bis}
$$
where $\rho(\cdot)$ identifies the the spectral radius of $R(\alpha)$. The underlying idea is that, in the vicinity of the convergence point $R(\alpha)$ becomes stationary, and therefore the most relevant parameter that determines the convergence speed is the spectral radius,\footnote{The role of the spectral radius is a consequence of the fact that $\rho(A) = \lim_{k\rightarrow\infty}\|A^k\|^{1/k}$ (e.g., see \cite{horn2013matrix}).} which satisfies $\rho(\cdot)\le\|\cdot\|$ and hence is a more precise indicator than the matrix norm. 

Note that \e{VAR52} is providing a \emph{global} information on the convergence speed, while \e{VAR52bis} is providing a \emph{local} information in the vicinity of the target point. Hence, both are valuable for understanding how to possibly control the convergence process. However, the spectral radius perspective is particularly welcome thanks to the fact that we are able to effectively locate the eigenvalues of $R(\alpha)$ starting from a loose indication given by the range constraint $\alpha_i\in[{\ell}_i, {\nu}_i]$. We show this by first concentrating on a simpler scenario where $\alpha_i$ can take two values only, which captures the essence of the problem, and by then generalising the result to a more conventional setting.

\subsection{A simplified scenario}\label{sec-locating}

We begin by investigating the span of the eigenvalues and the singular values of $R(\alpha)$ in a reference scenario where each $\alpha_i$ can take either a positive or a negative value (one of which can be zero). This scenario is, e.g., active at convergence for a linear problem, in which case the limit point $z^*$ lies somewhere on a staircase function of the form of \fig{CO4}.a and assumes values $1$ or $-1$. The more general framework that we propose, however, is one where levels can be different from $\pm1$, which generalises the results available in \cite{liang2017local} (see also \cite{bauschke2016optimal,demanet2016eventual,bjorck1973numerical}). 

The reference tool for obtaining our results is the CS decomposition of \cite[Theorem~2.5.3]{golub2012matrix} (see details in Appendix~\ref{App-teoremaAutovalFixed}), and provides the following

\begin{theorem}\label{teoremaAutovalFixed}
Assume that the entries of $\alpha$ in \e{VAR48} are ordered and take either a positive or a negative value, that is
$$
{\rm diag}(\alpha_i) =  H_i \triangleq  \qmatrix{ \overline{p}_i I_{p_i\times p_i}\cr & - \overline{n}_i I_{n_i\times n_i}}\;,
\e{MV10}
$$
for some positive constants $\overline{p}_i\ge0$ and $\overline{n}_i\ge0$, and where $p_i$ and $n_i$ denote, respectively, the number of positive and negative entries. Accordingly to the above, consider expressing the matrix product $V_1V_2^T$ in the form
$$
V_2V_1^T = G = \qmatrix{G_1 & G_2\cr G_3 & G_4}\;,
\e{MV12}
$$
where $G_1$ is $p_2\times p_1$, $G_2$ is $p_2\times n_1$, $G_3$ is $n_2\times p_1$, and $G_4$ is $n_2\times n_1$. Denote with $c_1,\ldots,c_p$ the singular values of either $G_1$ or $G_4$, where $G_1$ is used if  $\min(p_1,p_2)\le\min(n_1,n_2)$ and $G_4$ is used otherwise. The singular values are in number of $p = \min(p_1,p_2,n_1,n_2)$, and belong to the range $c_i\in[0,1]$. Define constants
$$
\eqalign{
\overline{n} & = \max(\overline{p}_1 \overline{n}_2, \overline{n}_1 \overline{p}_2)\;,\cr
\overline{p} & = \max(\overline{p}_1 \overline{p}_2, \overline{n}_1 \overline{n}_2)\;,\cr
\overline{r} & = \sqrt{\overline{p}_1 \overline{n}_1 \overline{n}_2\overline{p}_2}\;.\cr
\,
}\quad\eqalign{
\breve{n} & = \min(\overline{p}_1 \overline{n}_2, \overline{n}_1 \overline{p}_2)\;,\cr
\breve{p} & = \min(\overline{p}_1 \overline{p}_2, \overline{n}_1 \overline{n}_2)\;,\cr
k_1  &= \fract12  (\overline{p}_1 + \overline{n}_1) (\overline{p}_2 + \overline{n}_2)\;, \cr
k_2 & =  \fract12 (\overline{p}_1\overline{n}_2+ \overline{n}_1\overline{p}_2)\;.\cr
}
\e{MV11}
$$
Then the eigenvalues of matrix $N(\alpha)$ are known in the closed form, and we have that:
\begin{enumerate}
\item each singular value $c_i$ generates a couple of distinct eigenvalues by the rule
$$
\lambda_{1,2}(c_i) =  k_1 c_i^2 - k_2  \pm\sqrt{(k_1 c_i^2 - k_2)^2 - \overline{r}^2} \;,
\e{MV14}
$$
for a total of $2p$ eigenvalues;
\item the remaining eigenvalues are real--valued, and can be positive with value $\overline{p}$ or $\breve{p}$ (either of the two, but not both), or negative with value $-\overline{n}$ or $-\breve{n}$ (either of the two, but not both).
\end{enumerate}
\hfill~$\Box$
\end{theorem}
\begin{IEEEproof} See Appendix~\ref{App-teoremaAutovalFixed}.
\end{IEEEproof}

We observe that, depending on the value of $c_i$, the eigenvalues \e{MV14} can be real and positive, real and negative, or complex valued. The transition is controlled by the threshold values
$$
 \underline{c} =  \frac{|\sqrt{\overline{p}_1\overline{n}_2\rule{0mm}{3mm}}-\sqrt{\rule{0mm}{3mm}\overline{n}_1\overline{p}_2}\;|}{ \sqrt{2k_1}}\;,\quad
 \overline{c} = \frac{\sqrt{\overline{p}_1\overline{n}_2\rule{0mm}{3mm}}+\sqrt{\rule{0mm}{3mm}\overline{n}_1\overline{p}_2}\;}{ \sqrt{2k_1}}\;,
 \e{MV40}
$$
which satisfy $\underline{c}\le\overline{c}$, and which guarantee that the eigenvalues \e{MV14} satisfy the following:
\begin{enumerate}
\item[a)] for $c_i<\underline{c}$ they are real valued and \emph{negative} in the range $[-\overline{n},-\breve{n}]$;
\item[b)] for $c_i>\overline{c}$ they are real valued and \emph{positive} in the range $[\breve{p},\overline{p}]$;
\item[c)] for $c_i\in[\underline{c},\overline{c}]$ they form a \emph{complex--conjugate} couple with absolute value $\overline{r}$.
\end{enumerate}
According to this result, a pictorial representation of the locus of the eigenvalues of $N(\alpha)$, namely the set
$$
\M N = \Big\{{\rm eigs}(N(\alpha)), \alpha_{i,j}\in\{-\overline{n}_i,\overline{p_i}\} \Big\}\;,
\e{VF3}
$$
is shown in \fig{CO16}\Fig[t]{CO16}. 

Incidentally, the above characterisation is complete for an understanding of the spectral behaviour of matrix $R(\alpha)$, i.e., for identifying \e{VAR52bis}. As a matter of fact, the locus $\M N$ is mapped in the eigenvalue set of $R(\alpha)$ by a simple linear map, namely 
$$
\M R = \Big\{{\rm eigs}(R(\alpha)), \alpha_{i,j}\in\{-\overline{n}_i,\overline{p_i}\} \Big\} = (1-q) + q \M N\;,
\e{VF4}
$$
whose action is exemplified in \fig{CO9} for:
\begin{enumerate}
\item a convex problem (\fig{CO9}, above) where \e{VAR30} ensures that $\M N$ is contained in the unit circle, and showing how the map \e{VF4} can effectively reduce the spectral radius; we used $0<q<1$, which corresponds to a contraction in the complex plane towards the \emph{fixed} point $1$ (the point is named \emph{fixed} since the map $\M N\rightarrow\M R$ does not change its position, this being true for any choice of $q$);
\item a non--convex problem (\fig{CO9}, below) where the action of $q$ is unable to map inside the unit circle those points that, broadly speaking, lie on the right of the fixed point $1$.
\end{enumerate}
Note that the above clearly explains the role of $q$, that is, it controls the possibility of shifting the eigenvalues, possibly inside the unit circle and towards $0$ (in such a way to increase the convergence speed). It also suggests that different eigenvalues would need in principle a different value of $q$ since, e.g.,  $\lambda=0.95$ calls for a value $q>1$ to be shifted towards zero, while for a similar effect $\lambda=-0.95$ calls for $0<q<1$, and $\lambda=1.05$ calls for $q<-1$. \Fig[t]{CO9}

From a graphical inspection of \fig{CO9} we can also observe that, irrespective of the value of $q$, the spectral radius is determined by the values of $\M R$ lying on the real axis, the reasons being that \e{VF4} maps the circle of radius $\overline{r}$ in a circle of radius $|q|\overline{r}$ centred in $1-q$. This provides the following results whose (simple) proof is left to the reader:
\begin{corollary}\label{corollarioBestq}
Under the assumptions of Theorem~\ref{teoremaAutovalFixed}, the spectral radius \e{VAR52bis} is upper bounded by
$$
\rho \le \rho_{\max} = \max\Big( |(1-q)-q\overline{n}|,|(1-q)+q\overline{p}| \Big)\;.
\e{NEC9}
$$
The resulting optimum choice $q_\rho^*$ minimising \e{NEC9} is therefore
$$
q_\rho^* =\cases{\displaystyle\frac2{2+\overline{n}-\overline{p}}\;,\quad & $0<\overline{p}<1$\cr
\rule{0mm}{4mm}0\;,\quad & $\overline{p}\ge1$\cr
}
\e{NEC9b1}
$$
and the related bound is
$$
\rho_{\max}^* =\cases{\displaystyle \frac{\overline{p}+\overline{n}}{2+\overline{n}-\overline{p}}\;,\quad & $0<\overline{p}<1$\cr
\rule{0mm}{4mm}1\;,\quad & $\overline{p}\ge1$.\cr
}
\e{NEC9b}
$$
\hfill~$\Box$
\end{corollary}
\begin{IEEEproof} Left to the reader.
\end{IEEEproof}
Interestingly, for $0<\overline{p}<1$ the bound \e{NEC9b} satisfies $\rho_{\max}^*< 1$, in which case we are guaranteeing a locally linear convergence to the target point. When, instead, $\overline{p}>1$, then the above formalisation is suggesting the inappropriate choice $q=0$, that is, we are not able to build a convergent algorithm.

\subsection{A more general statement}\label{sec-general}

Although the four--level assumption of \e{MV10} is linked to a very specific setting observed at convergence, its outcomes are more general than one can expect. As a matter of fact, the structures identified so far (i.e., the shape of \fig{CO16}) are observed also when the values of $\alpha$ belong to a range. This leads to our general statement, providing a description of the convergence process with a larger view. 

\begin{theorem}\label{teoremaAutoval}
Assume that the entries of $\alpha_i$ in \e{VAR48} are limited to the range
$$
\alpha_{i,j} \in [-\overline{n}_i,-\underline{n}_i] \cup [\underline{p}_i,\overline{p}_i]\;, \quad \forall j
\e{MV2}
$$
with $0\le\underline{n}_i\le\overline{n}_i$ and $0\le\underline{p}_i\le\overline{p}_i$, and define the positive constants 
$$
\eqalign{
\overline{n} & = \max(\overline{p}_1\overline{n}_2,\overline{n}_1\overline{p}_2)\;,\cr
\overline{p} & = \max(\overline{p}_1\overline{p}_2,\overline{n}_1\overline{n}_2)\;,\cr
\overline{r} & = \sqrt{\overline{p}_1\overline{p}_2\overline{n}_1\overline{n}_2}\;,
}\quad \eqalign{
\underline{n} & = \min(\underline{p}_1\underline{n}_2,\underline{n}_1\underline{p}_2)\;,\cr
\underline{p} & = \min(\underline{p}_1\underline{p}_2,\underline{n}_1\underline{n}_2)\;,\cr
\underline{r} & = \sqrt{\underline{p}_1\underline{p}_2\underline{n}_1\underline{n}_2}\;,
}
\e{MV4}
$$
with $0\le\underline{n}\le \underline{r} \le \overline{r} \le \overline{n}$ and $0\le\underline{p}\le \underline{r} \le \overline{r} \le \overline{p}$. Then, the real--valued eigenvalues of $N(\alpha)$ belong to the range
$$
\lambda \in [-\overline{n},-\underline{n}] \cup [\underline{p},\overline{p}]\;,
\e{MV6}
$$
while the complex--valued eigenvalues satisfy
$$
\underline{r} \le |\lambda| \le \overline{r}\;.
\e{MV8}
$$
\hfill~$\Box$
\end{theorem}
\begin{IEEEproof} See Appendix~\ref{App-teoremaAutoval}.
\end{IEEEproof}

The outcomes of Theorem~\ref{teoremaAutoval} are illustrated in \fig{CO17}\Fig[t]{CO17}, where
$$
\M N = \Big\{{\rm eigs}(N(\alpha)), \alpha_{i,j}\in[-\overline{n}_i,-\underline{n}_i] \cup [\underline{p}_i,\overline{p}_i] \Big\}\;.
\e{VF3b}
$$
The implications on $R(\alpha)$ are, evidently, identical to the simplified case, that is, both \e{VF4} and Corollary~\ref{corollarioBestq} hold also in the general case. Notably, the structure of \fig{CO17} is recurring in many fields of study (e.g., in network science \cite{georgeot2010spectral} or in the characterisation of large random matrices), 
the reason being that with Theorem~\ref{teoremaAutoval} we are providing a characterisation of the eigenvalue span of any real matrix since any real matrix can be expressed as the product of two symmetric matrices \cite{horn2013matrix}.  

The interesting outcome of our result is that by identifying the parameters $\overline{n}_i,\underline{n}_i,\underline{p}_i,\overline{p}_i$ in both the alternate directions (whose identification usually corresponds to a very simple task), we are able to characterise the location of the eigenvalues of $R(\alpha)$ without the need to calculate them explicitly, the added value being that calculating the eigenvalues of a large matrix might be a costly operation.

\Section[EX]{Application example}

For illustration purposes we consider the simple weighted Lasso minimisation problem used in \cite{giselsson2017linear}
$$
\min_x \;\underbrace{\fract12\|\Omega \,x - o\|_2^2}_{f_1} + \underbrace{\|w\circ x\|_1}_{f_2}
\e{LA2}
$$
where $\|\cdot\|_p$ is the Euclidean norm of order $p$, $\circ$ denotes the entry--wise Hadamard product, $x$, $o$, and $w$ are real--valued vectors of length, respectively, $200$, $300$, and $200$, where $\Omega$ is a sparse $300\times200$ data matrix with an average of $10$ active entries per row, and where the (active) entries of $\Omega$ and $\omega$ are drawn from a Gaussian distribution with zero mean and unit variance, while the entries of $w$ are taken from a uniform distribution in $[0,1]$. 

The problem is solved through \ac{ADMM} by setting $x_1=x_2=x$, hence $A_1=A_2=I$ and $b=0$. We also assume a diagonal augmentation matrix $E=\epsilon I$. This provides a linear reflected proximity operator $\C D_1$ of the form
$$
\C D_1(u) =\left( \frac{I-\Omega^T\Omega \,\epsilon^{-1}}{I+\Omega^T\Omega\,\epsilon^{-1}} \right)\cdot u + \frac{2\Omega^To/\sqrt{\epsilon}}{I+\Omega^T\Omega\,\epsilon^{-1}}\;,
\e{LA6}
$$
where $\Omega^T\Omega\,\epsilon^{-1}=\tilde{S}_1=\tilde{C}_1$, so that in \e{HH10} we have $\ell_{1,k}=\nu_{1,k}=h(\lambda_k(\Omega^T \Omega)/\epsilon)$. This corresponds to a range of the form $\alpha_{1,j}  \in [h(\lambda_{\max}/\epsilon),h(\lambda_{\min}/\epsilon)]$ where
$$
\lambda_{\max}  = \max_k \lambda_k(\Omega^T \Omega)\;,\qquad \lambda_{\min}  = \min_k \lambda_k(\Omega^T \Omega)\;,
\e{LA8}
$$
are the maximum and minimum eigenvalue, respectively, of matrix $\Omega^T \Omega$. As a straightforward consequence, we can set
$$
\eqalign{
\overline{p}_1 & = [h(\lambda_{\min}/\epsilon)]^+\;,\qquad
\overline{n}_1  = [-h(\lambda_{\max}/\epsilon)]^+ \cr
\underline{p}_1 & = [h(\lambda_{\max}/\epsilon)]^+\;,\qquad
\underline{n}_1  = [-h(\lambda_{\min}/\epsilon)]^+ \;,
}
\e{LA8b}
$$
where $[u]^+$ is the plus operator providing $u$ for $u>0$ and $0$ otherwise. 
%
%
For the reflected proximity operator $\C D_2$ we can instead exploit the technique developed in \fig{CO4}, to derive its $j$th entry $\C D_{2,j}(u_j)$ as illustrated in \fig{CO24}\Fig[t]{CO24}. The staircase structure suggests that
$$
\overline{p}_2=\overline{n}_2 = 1\;,\qquad \underline{p}_2=\underline{n}_2 = 0\;,
\e{LA8c}
$$
on the entire iteration span, but $\underline{p}_2$ and $\underline{n}_2$ are ensured to be $\simeq 1$ in the vicinity of the target point. By using the weaker (but simpler) constraint \e{LA8c}, we obtain 
$$
\overline{r} = \sqrt{\overline{p}_1\overline{n}_1}\;,\qquad
\overline{p} = \overline{n}  =\max(\overline{p}_1,\overline{n}_1)
\e{LA10}
$$
and $\underline{p}=\underline{n}=\underline{r}  = 0$, whose symmetry, according to \e{NEC9b1}, identifies the optimum parameter $q=1$. 

The study of the algorithm convergence when using $q=1$ and $\epsilon=1$ is shown in \fig{CO28}\Fig[t]{CO28} and \fig{CO26}\Fig[t]{CO26} for one specific realisation of the random matrices and vectors $\Omega$, $o$, and $w$. In the specific example we have $\lambda_{\max}=60.75$ and $\lambda_{\min}=0.3465$, providing $\overline{p}_1= 0.4853$ and $\overline{n}_1= 0.9676$. Note in \fig{CO28} the correspondence between the locus of eigenvalues $\M R$ (in gray) and the eigenvalues at the limit point $z^*$ (dots). The locus $\M R$ was identified according to Theorem~\ref{teoremaAutoval} and \e{LA10}, and perfectly captures the position of the eigenvalues, especially in the real axis where the eigenvalue with the largest absolute value $\lambda^*_{\max}=0.9543$ is very close to the maximum allowed value $\overline{p}=0.9676$. Note that the bound given by $\M R$ is in any case much stricter in the real axis than in the complex plane (circle $\overline{r}$), this being a common behaviour for randomly wired networks such as the one we are considering. The convergence plot of \fig{CO28} further evidences that, although $\lambda^*_{\max}$ provides a precise indication of the convergence speed at the limit point, the value of $\rho_{\max}$ identified according to Theorem~\ref{teoremaAutoval} and \e{LA10}, namely $\rho_{\max} =0.9676$ in the considered case, is an upper value that is often able to capture with greater accuracy the true speed of the convergence process. In addition, $\rho_{\max}$ is easier to calculate since it does not rely on the (precise) knowledge of the limit point $z^*$, which in a practical application setup might be hardly known or derived.

\Section[CO]{Conclusions}

In this paper we provided new theoretical insights on the convergence process of \ac{ADMM}, which were derived by interpreting \ac{ADMM} as a (generalised) Douglas-Rachford splitting in the time domain (as opposed to the classical frequency domain interpretation). The interesting outcome is that convergence speed is linked to the eigenvalues of matrices $R(\alpha)$ whose general structure is evidenced in \fig{CO17}, where a constraint on the absolute value of complex eigenvalues cohabits with a (looser) constraint on the range of real-valued eigenvalues. The latter constraint is what will ultimately set the convergence speed, and calls for novel speeding methods. For example, under the assumption that $\partial \C D_i$ is constant, one might consider a cyclostationary $q$ where different values $q_k$ (cyclically repeating) are used at each iteration in order to selectively trim the effect of the larger (real-valued) eigenvalues. The solution in a more general and realistic setting is a tough problem which calls for future investigations and novel ideas.

\appendices

\section{Proof of Theorem~\ref{teoremaAutovalFixed}}\label{App-teoremaAutovalFixed}

In the following we will denote with $A \sim B$ a similarity between matrices $A$ and $B$ in the sense that they share the same eigenvalues with the same multiplicities. We will also use the property $AB\sim BA$. The starting point is to identify the similarity (see \e{VAR48})
$$
N(\alpha) \sim  \underbrace{V_2V_1^T}_{G} \cdot H_1 \cdot \underbrace{V_1V_2^T}_{G^T}\cdot H_2 \;,
\e{KH2}
$$
where we used the notation \e{MV10} and \e{MV12}. Note that matrices $H_i$ are real valued and diagonal, and matrix $G$ is \emph{orthogonal} (real valued and unitary). 

We then want to identify an operative structure for matrix $G$. To do so we initially assume that 
$$
p = p_2\le p_1\le n_1\le n_2\;.
\e{MB2}
$$ 
The general case, which is a simple extension of the outcomes of \e{MB2}, will be discussed at the end of the proof. Hence, we can readily exploit the CS decomposition of \cite[Theorem~2.5.3]{golub2012matrix}, which we write in the following form:

\begin{proposition}\label{golubMat}
Under condition \e{MB2}, matrix $G$ can be expressed as the product
$$
G  = \underbrace{\qmatrix{A_1\cr &A_2}}_{A} \cdot M \cdot \underbrace{\qmatrix{B_1\cr &B_2}}_{B}\;,
\e{MA20}
$$
where $A_1$, $A_2$, $B_1$, and $B_2$ are orthogonal matrices of dimension, respectively, $p_2\times p_2$, $n_2\times n_2$, $p_1\times p_1$, $n_1\times n_1$, and where $M$ is the orthogonal matrix 
$$
M  = \qmatrix{
\begin{picture}(0,0)(1,3)\put(0,0){\begin{tikzpicture}\draw (0,0) -- ++(4,0)[dotted];\end{tikzpicture}}\end{picture}
C & &
\begin{picture}(0,0)(1,40)\put(0,0){\begin{tikzpicture}\draw (0,0) -- ++(0,-1.7)[dotted];\end{tikzpicture}}\end{picture}
 -S\cr
S & & C \cr 
&  I_{q_1\times q_1} & & 0_{q_1\times q_2}\cr 
& 0_{q_2\times q_1} & & I_{q_2\times q_2}}
\e{MA20b}
$$
where $q_1 = p_1-p_2$ and $q_2=n_2-p_1=n_1-p_2$, and where $C = {\rm diag}(c_1,\ldots,c_p)$ and $S=\sqrt{I-C^2}$ are $p_1\times p_1$ diagonal, positive definite, and real valued matrices with entries in the range $[0,1]$. \hfill~$\Box$
\end{proposition}

Use of \e{MA20} in \e{KH2} provides
$$
\eqalign{
N(\alpha)  & \sim M B H_1 B^T \cdot M^T \cdot A^T H_2 A  \cr
 & = M B B^T \cdot H_1  M^T H_2 \cdot A^T  A  \cr
 & = \underbrace{M H_1  M^T H_2}_H   \cr
}
\e{MB4}
$$
where we exploited the fact that $H_1$ commutes with $B^T$, and $H_2$ commutes with $A^T$, as well as the fact that both $A$ and $B$ are orthogonal matrices. By substituting \e{MA20b}, we can then explicitly express matrix $H$ in \e{MB4} in the block--diagonal form
$$
H = \qmatrix{
H_0 \cr
& -\overline{p}_1\overline{n}_2 I_{q_1\times q_1} \cr
& & \overline{n}_1\overline{n}_2 I_{q_2\times q_2}
}\;,
\e{MB10}
$$
where
$$
H_0 = \qmatrix{
(C^2 \overline{p}_1 - S^2\overline{n}_1)\overline{p}_2 &
-CS (\overline{p}_1 + \overline{n}_1)\overline{n}_2\cr
CS (\overline{p}_1 + \overline{n}_1)\overline{p}_2 & 
(C^2 \overline{n}_1 - S^2\overline{p}_1)\overline{n}_2 \cr
}\;.
\e{MB12}
$$
The eigenvalues in \e{MV14} are exactly the roots of $H_0$, which can be identified by standard techniques since \e{MB12} has a simple $2\times 2$ block structure.\footnote{We skip the derivation since it is trivial, but tedious, and only report the final result.} By inspection of \e{MB10}, the remaining eigenvalues are $\overline{n}_1\overline{n}_2$ and $-\overline{p}_1\overline{n}_2$, which is consistent with statement 2) given by the theorem.

The extension to the cases where \e{MB2} is not valid can be performed with a simple substitution of variables. To start with, consider the case where $\min(\overline{p}_2,\overline{n}_2)\le \min(\overline{p}_1,\overline{n}_1)$, which entails four different scenarios, namely:
\begin{enumerate}
\item[a)] $p = p_2\le p_1\le n_1\le n_2$, which we already discussed;
\item[b)] $p = p_2\le n_1\le p_1\le n_2$,
\item[c)] $p = n_2\le p_1\le n_1\le p_2$, and
\item[d)] $p = n_2\le n_1\le p_1\le p_2$.
\end{enumerate}

Now, scenario b) can be solved equivalently to scenario a) provided that the matrix $G$ is rearranged in the form
$$
\qmatrix{G_2 & G_1\cr G_4 & G_3}\;.
\e{MBR2}
$$
This corresponds to making the swaps $p_1\leftrightarrow n_1$ and $\overline{p}_1\leftrightarrow \overline{n}_1$, and to multiplying by $-1$ the resulting eigenvalues, since the role of positive and negative have been swapped in $H_1$. Incidentally, \e{MBR2} implies that the singular values of interest are now those of $G_2$. Thanks to the orthogonality of $G$ we have $G_2G_2^*=I-G_1G_1^*$, and therefore the squared singular values of $G_2$ are $1-c_i^2$. This implies a further substitution $c_i^2\rightarrow 1-c_i^2$, i.e., the swap $c_i\leftrightarrow s_i$. Notably, under the swap $\overline{p}_1\leftrightarrow \overline{n}_1$, the swap $c_i\leftrightarrow -s_i$, and the multiplication by $-1$, the expression \e{MV14} remains valid. As a matter of fact, we are simply transforming $H_0$ into $H_0^T$. The other eigenvalues of \e{MB10}, accordingly, become $\overline{n}_1\overline{n}_2$ and $- \overline{p}_1\overline{n}_2$, which is consistent with statement 2).

Scenarios c) and d) can be solved equivalently. 
%
%
%
The remaining cases to consider are those where $\min(\overline{p}_2,\overline{n}_2)\ge \min(\overline{p}_1,\overline{n}_1)$, in which case we simply need to swap the roles of $i=1$ and $i=2$, and to work with $G^T$, but this has no consequences on \e{MV14}, nor on the statement on the remaining eigenvalues.

\section{Proof of Theorem~\ref{teoremaAutoval}}\label{App-teoremaAutoval}

In the following, the result \e{MV8} on the range of complex eigenvalues, and the result \e{MV6} on the range of real eigenvalues, are separately proved since they require different techniques. Following Appendix~\ref{App-teoremaAutovalFixed}, we denote with $A \sim B$ a similarity between matrices $A$ and $B$ in the sense that they share the same eigenvalues with the same multiplicities, and extensively use the property $AB\sim BA$. The starting point is to identify the similarity (see \e{VAR48})
$$
N(\alpha) \sim  G  \Lambda_1 G ^T \Lambda_2 \;,
\e{KH2bis}
$$
where $G=V_2V_1^T$ is as orthogonal matrix, and where $\Lambda_i={\rm diag}(\alpha_i)$ with $\alpha_i$ in the ranges \e{MV2}. We also prove the theorem by considering $\underline{n}_i>0$ and $\underline{p}_i>0$, i.e., that matrices $\Lambda_i$ are invertible. The applicability to the case where $\underline{n}_i=0$ and/or $\underline{p}_i=0$ is guaranteed by taking a limit thanks to the results of perturbation theory (e.g., see \cite{horn2013matrix}).

\subsection{Complex eigenvalues -- Proof of \e{MV8}}

Let $x$ be a complex--valued eigenvalue of matrix $G \Lambda_1 G^T \Lambda_2$, with unit norm $\|x\|=1$, and with associated eigenvector $\lambda$. This allows writing the general relation
$$
\Lambda_2 x = \lambda \cdot G \Lambda_1^{-1} G^T x\;,
\e{TG2}
$$
which is true thanks to the invertibility of $\Lambda_1$. In turn, from \e{TG2} we further obtain
$$
\eqalign{
x^*\Lambda_2 x & = \lambda \cdot x^*G \Lambda_1^{-1} G^T x\cr
x^*\Lambda_2^2 x & = |\lambda|^2 \cdot x^*G \Lambda_1^{-2} G^T x\;,
}
\e{TG4}
$$
the first being derived by multiplying \e{TG2} by $x^*$, the second by taking the squared norm. From the second of \e{TG4} we have
$$
|\lambda|^2 = \frac{x^*\Lambda_2^2 x}{x^*G \Lambda_1^{-2} G^T x} > 0\;,
\e{TG6}
$$
which is ensured by the invertibility of $\Lambda_1$ and $\Lambda_2$ (it implies that  both $\Lambda_2^2$ and $\Lambda_1^{-2} $ are positive definite, so that both the numerator and the denominator are strictly positive thanks to the assumption $\|x\|=1$). A perfectly equivalent outcome cannot be stated for the first of \e{TG4}, since $x^*G \Lambda_1^{-1} G^T x$ can be equal to zero. In this case we instead have two distinct possibilities, namely: 1)
$$
\lambda = \frac{x^*\Lambda_2 x}{x^*G \Lambda_1^{-1} G^T x}\quad\hbox{subject to } x^*G \Lambda_1^{-1} G^T x\neq 0 \;,
\e{TG8}
$$
which clearly identifies real--valued eigenvalues, and 2)
$$
\left\{\eqalign{
&x^*\Lambda_2 x = 0 \cr
&x^*G \Lambda_1^{-1} G^T x = 0\;,
}\right.
\e{TG10}
$$
which in general identifies complex--valued eigenvalues. 

A range for the complex valued eigenvalues can therefore be identified by considering both \e{TG6} and \e{TG10}, namely the optimisation problem
$$
\eqalign{
&\maxim \frac{x^*\Lambda_2^2 x}{y^* \Lambda_1^{-2} y}\cr
& \st \|x\|=1\;,\; y= G^T x\;,\; \|y\|=1\;,\cr
& \pst x^*\Lambda_2 x = 0\;,\; y^* \Lambda_1^{-1} y = 0\;,
}
\e{TG12}
$$
which clearly provides and upper bound to $|\lambda| ^2$. As a simpler alternative we consider relaxing \e{TG12} by removing the constraint $y= G^T x$, which provides two simple optimisation problems of the form
$$
\eqalign{
&\maxim x^*\Lambda_2^2 x\quad \st \|x\|=1\;,\;  x^*\Lambda_2 x = 0\;,
}
\e{TG14}
$$
and
$$
\eqalign{
&\minim y^* \Lambda_1^{-2} y\quad \st \|y\|=1\;,\;  y^* \Lambda_1^{-1} y = 0\;,
}
\e{TG16}
$$
whose composition still ensures an upper bound to $|\lambda| ^2$. 

We preliminarily concentrate on the solution to \e{TG14}, which we do by denoting $|x_i|^2=a_i$, so that the problem can be equivalently written in the form
$$
\eqalign{
&\maxim \sum_i a_i \lambda_{2i}^2 \cr
& \st a_i\ge0\;,\; \sum_i a_i \lambda_{2i} = 0\;,\; \sum_ia_i=1\;.
}
\e{TG18}
$$
The corresponding Lagrangian assumes the form
$$
L(a,\mu) = \sum_i \Big[a_i \big(\lambda_{2i}^2  + \mu_1 \lambda_{2i} + \mu_2 \big) + \eta_{\ge 0}(a_i)\Big]- \mu_2\;,
$$
with $\eta_{\ge 0}(x)$ an indicating function for $x\ge0$. Observe that the \ac{KKT} condition $a^* = \max_a L(a,\mu^*)$, where $(a^*, \mu^*)$ denotes the primal--dual optimal points, ensures that 
$$
g^*_i = \lambda_{2i}^2  + \mu_1^* \lambda_{2i} + \mu_2^* \le 0\;,\quad \forall i\;.
\e{TG20}
$$ 
In fact, $g^*_i<0$ implies $a_i^*=0$, $g^*_i=0$ implies any value for $a_i^*$, while $g^*_i>0$ implies $a_i^*=+\infty$ which leads to an infeasible solution. Actually, the \ac{KKT} condition $\sum_i a_i^*=1$ implies that there must be at least one $g^*_i = 0$, while the \ac{KKT} condition $\sum_i a_i^* \lambda_{2i}=0$ further implies that the active values must be at least two, one for a positive eigenvalue $ \lambda_{2i}$, and one for a negative eigenvalue. Having said so, it is straightforward to see that the quadratic form of $g^*_i$ as a function of $\lambda_{2i}$ is required to have roots in $\overline{p}_2$ and $-\overline{n}_2$, in such a way that $g_i^*=0$ for those $i$ where $\lambda_{2i}=\overline{p}_2$ or $-\overline{n}_2$, and $g_i^*<0$ otherwise. As a consequence, it is $a_i^* = 0 \, \forall i \hbox{ such that } \lambda_{2i}\neq \overline{p}_2,-\overline{n}_2$. By then using the compact notation
$$
\alpha = \sum_{i|\lambda_{2i}= \overline{p}_2} a_i^* = 1- \sum_{i|\lambda_{2i}= -\overline{n}_2} a_i^*
\e{TG22}
$$
we can reinterpret \e{TG18} in the form
$$
\eqalign{
&\maxim \alpha \overline{p}_2^2 + (1-\alpha) \overline{n}_2^2\cr
&\st 0\le\alpha\le 1, \;\alpha \overline{p}_2 = (1-\alpha) \overline{n}_2 
}
\e{TG24}
$$
which can be easily solved by exploiting the second constraint, that is by replacing $\alpha = (1-\alpha) \overline{n}_2/\overline{p}_2$ and $(1-\alpha) = \alpha  \overline{p}_2/\overline{n}_2$ in the target function. The upper bound in \e{TG12} is therefore simply $\overline{p}_2\overline{n}_2$.

The solution to \e{TG14} is perfectly equivalent, with the observation that the counterpart to \e{TG20} is 
$$
g^*_i = \lambda_{1i}^{-2}  + \mu_1^* \lambda_{1i}^{-1} + \mu_2^* \ge 0\;,\quad \forall i\;,
\e{TG20b}
$$ 
and leads to a choice of $\mu^*$ that sets $g^*_i=0$ for those $i$ associated with eigenvalues $\overline{p}_1^{-1}$ and $-\overline{n}_1^{-1}$, namely the least positive and greatest negative values, and $g_i^*>0$ for the rest. This identifies a lower bound $1/(\overline{p}_1\overline{n}_1)$ in \e{TG14}.

By combining the two bounds we therefore have $|\lambda|^2 \le \overline{p}_1\overline{n}_1\overline{p}_2\overline{n}_2$, which is the upper bound in \e{MV6}. The lower bound can be stated equivalently, by simply investigating the upper bound of the inverse matrix 
$$
(G \Lambda_1 G^T \Lambda_2)^{-1}= \Lambda_2^{-1} G \Lambda_1^{-1} G^T \sim G \Lambda_1^{-1} G^T \Lambda_2^{-1}\;,
\e{TG30}
$$
where the roles of $\overline{p}_i$ and $\overline{n}_i$ swap with $\underline{p}_i^{-1}$ and $\underline{n}_i^{-1}$.

\subsection{Real eigenvalues -- Proof of \e{MV6}}

For the real eigenvalues case we exploit the results readily available from Appendix~\ref{App-teoremaAutovalFixed}, and assume that matrices $\Lambda_i$ are organised in such a way that the positive values are in the first positions, and are followed by the negative values. This allows writing $\Lambda_i$ in the form 
$$
\Lambda_i = \breve{\Lambda}_i H_i\;,
\e{MB7}
$$
where $H_i$ was defined in \e{MV10} and carries the largest positive and negative values of $\Lambda_i$, and where $\breve{\Lambda}_i$ is a diagonal positive definite matrix satisfying $0\prec \breve{\Lambda}_i\preceq I$. By replacing the equivalence in \e{KH2bis}, and by exploiting the commutative property between some of the considered matrices, we find that
$$
\eqalign{
N(\alpha)  & \sim M B \breve{\Lambda}_1 B^T \cdot H_1 M^T H_2 \cdot A^T \breve{\Lambda}_2 A \cr
& \sim \underbrace{M B \breve{\Lambda}_1 B^T M^T }_{P_1} 
	\cdot \underbrace{M H_1 M^T H_2 }_H \cdot \underbrace{A^T \breve{\Lambda}_2 A}_{P_2}\cr
}
\e{MB8}
$$
where $P_i$ are by construction real--valued positive--definite symmetric matrices with $0\prec P_i\preceq I$, and where $H$ is the matrix defined in \e{MB10}-\e{MB12}. 

Now, from Appendix~\ref{App-teoremaAutovalFixed} we know that $H$ is real and diagonalizable, with eigenvalues belonging to the range depicted in \fig{CO16} (see details in Section~\ref{sec-locating}). This ensures that it can be written in the form\footnote{$H$ is diagonalisable, hence it can be expressed in the form $U D U^{-1}$ with diagonal $D$ collecting the eigenvalues. The columns of $U$ are the eigenvectors, that can be chosen real--valued for real--valued eigenvalues. For complex valued eigenvalue couples $\{\lambda,\lambda^*\}$, instead, the eigenvectors can be chosen to be a complex conjugate couple, i.e., $\{u,u^*\}$. To obtain a real--valued form it suffices multiplying the blocks of interest by a normalising unitary matrix $J$, namely
$$
\underbrace{[u, u^*]J }_{V}  \cdot \underbrace{J^*\qmatrix{\lambda\cr & \lambda^*} J}_B \cdot J^* \;,\quad J = \frac1{\sqrt{2}} \qmatrix{1 & i\cr 1& -i}
$$  
where matrix $V$ is now real--valued, and where matrix $B$ assumes the real--valued form \e{MB23}.
}
$$
H = V_H \Lambda_H V_H^{-1}\;,
\e{MB22}
$$
where $V_H$ is real--valued and invertible, while $\Lambda_H$ is a real--valued and normal matrix with a block diagonal structure whose blocks are either:
\begin{enumerate}
\item $1\times 1$ blocks carrying the real--valued eigenvalues of $H$;
\item $2\times 2$ blocks of the form
$$
B(c_i) = \qmatrix{\Re\lambda_{1}(c_i) & - \Im \lambda_{1}(c_i)\cr\Im \lambda_{1}(c_i)& \Re \lambda_{1}(c_i)}
\e{MB23}
$$
for all the complex--conjugate couples in \e{MV14}. 
\end{enumerate}
Hence, we can write
$$
N(\alpha)  \sim \underbrace{V_H^{-1} P_2 P_1 V_H}_T \cdot \Lambda_H  \;,
\e{MB24}
$$ 
which reveals the fundamental role of matrix $T$.

Notably, matrix $T$ in \e{MB24} is a \emph{weakly positive} matrix, that is, a diagonalisable matrix with non-negative eigenvalues, a result which is a straightforward consequence of the fact that the product of two positive definite matrices is weakly positive \cite{wigner1993weakly,nilssen2005weakly}. Being weakly positive, $T$ satisfies
$$
x^T T x\ge0 \;,\qquad \forall x\in\M R^n\;.
\e{MB32a}
$$
Incidentally, the bounds $P_1\preceq I$ and $P_2\preceq I$ guarantee that $\|P_2P_1\|\le 1$. Since $T\sim P_2P_1$, the above further implies that the largest eigenvalue of $T$ is bounded by $\lambda_{\max}(T)\le 1$. Similarly, the invertibility of both $P_1$ and $P_2$ further guarantee the existence of matrix $T^{-1}$, as well as the fact that its smaller eigenvalue satisfies $\lambda_{\min}(T^{-1})\ge 1$. This implies the property
$$
x^T T^{-1} x\ge \|x\|^2 \;,\qquad \forall x\in\M R^n\;,
\e{MB32b}
$$
which can be simply derived by observing that we can express $T^{-1}$ as $T^{-1}=I+T'$ with $T'$ some weakly positive matrix.

Given the above, a range for real eigenvalues can be derived by investigating the relation
$$
\Lambda_H x = \lambda T^{-1} x\;, \qquad \|x\|=1
\e{MB40}
$$
where $x$ is a real--valued eigenvector of unit norm, and $\lambda$ is its corresponding real eigenvalue. By multiplying by $x^T$, \e{MB40} allows writing 
$$
\lambda = \frac{x^T\Lambda_H x}{x^T  T^{-1} x}\;,
\e{MB42}
$$
where the denominator is ensured by \e{MB32b} to satisfy $x^T  T^{-1} x \ge 1$, while for the numerator inspection of \fig{CO16} ensures $x^T\Lambda_H x\in[-\overline{n},\overline{p}]$. Overall, \e{MB42} justifies the extrema of \e{MV6}, namely the relation $\lambda \in [-\overline{n},\overline{p}]$.

The relation $\lambda^{-1} \in [-\underline{n}^{-1},\underline{p}^{-1}]$, which completes the proof, can be derived by an equivalent method that investigates the inverse matrix \e{TG30}.

\bibliographystyle{IEEEtran}
\bibliography{co}

\end{document}